\documentclass{article}

\usepackage{amsfonts}
\usepackage{amsmath}
\usepackage{amssymb}
\usepackage{amsthm}
\usepackage{amscd}
\usepackage{verbatim}
\usepackage[dvips]{graphicx}
\usepackage{graphics}
\usepackage[all,2cell]{xypic}

\theoremstyle{plain}
\newtheorem{thm}{Theorem}[section]

\newtheorem{cor}[thm]{Corollary}
\newtheorem{lem}[thm]{Lemma}
\newtheorem{prop}[thm]{Proposition}

\theoremstyle{definition}

\theoremstyle{remark}

\DeclareMathOperator{\isoto}{\overset{\scriptstyle{\sim}}{\to}}

\DeclareMathOperator{\coker}{Coker}

\DeclareMathOperator{\Codim}{Codim}
\DeclareMathOperator{\Spec}{Spec}
\DeclareMathOperator{\Supp}{Supp}
\DeclareMathOperator{\Perf}{Perf}

\DeclareMathOperator{\LAX}{\mathcal{LAX}}
\DeclareMathOperator{\Aut}{Aut}

\DeclareMathOperator{\colim}{colim}

\DeclareMathOperator{\Homo}{H}

\DeclareMathOperator{\FFF}{\mathcal{F}}

\DeclareMathOperator{\KKK}{\mathcal{K}}
\DeclareMathOperator{\MMM}{\mathcal{M}}
\DeclareMathOperator{\PPP}{\mathcal{P}}

\DeclareMathOperator{\ulG}{\underline{G}}

\DeclareMathOperator{\Sh}{Sh}

\DeclareMathOperator{\lax}{lax}

\DeclareMathOperator{\strict}{s}

\DeclareMathOperator{\edge}{e}

\DeclareMathOperator{\ord}{ord}

\UseAllTwocells

\title{Gersten's conjecture}
\author{Satoshi Mochizuki}
\date{}

\begin{document}

\maketitle

\begin{abstract}
The purpose of this article is to prove that Gersten's conjecture for a 
commutative regular local ring is true. As its applications, we will prove the vanishing conjecture for certain Chow groups, generator conjecture for certain $K$-groups and Bloch's formula for absolute case.
\end{abstract}

\tableofcontents

\section{Introduction}

\noindent
The purpose of this note is to prove the following theorem.

\begin{thm}[Gersten's conjecture] \label{Main theorem}  $ $
\\
For any commutative regular local ring $R$, Gersten's conjecture is true.That is for any natural numbers $n$, $p$, the canonical inclusion 
$\MMM^{p+1}(R) \hookrightarrow \MMM^p(R)$ induces the zero 
map on $K$-groups
$$K_n(\MMM^{p+1}(R)) \to K_n(\MMM^p(R)) \\\ ,$$
where $\MMM^i(R)$ is the category of finitely generated $R$-modules $M$
with\\
 $\Codim_{\Spec R} \Supp M \geqq i$. 
\end{thm}

\noindent
Gersten's conjecture is proposed in \cite{Ger73}. More precise historical back grounds of this conjecture are explained in \cite{Moc07}. In \S 1, we will prove the main theorem and in \S 2, we will also discuss applications of this conjecture.

\noindent
\textbf{Acknowledgement} The author thankful to Shuji Saito for encouraging him, to Fabrice Orgogozo for stimulating argument about Corollary \ref{gen con}, to Takeshi Saito for making him to get to the reduction argument in Lemma \ref{Ret pri}, and to Kazuhiko Kurano for teaching him condition (iii) in Proposition \ref{equi cond}.

\section{Proof of the main theorem}

\noindent
From now on, let $R$ be a commutative regular local ring. Proof of the main theorem is divided series of lemmas. First we will improve Quillen's reduction argument in the proof of Gersten's conjecture in \cite{Qui73}.

\begin{lem}[Quillen induction] \label{Qui ind} $ $
\\
To prove the main theorem, we shall only check the following assertion:\\
For any non-negative integers $n$, $p$, and any regular sequence $f_1,\ldots,f_{p+1}$ in $R$, the canonical map induced from the inclusion map $\PPP(R/(f_1,\ldots,f_{p+1})) \hookrightarrow \MMM^p(R)$,
$$K_n(\PPP(R/(f_1,\ldots,f_{p+1}))) \to K_n(\MMM^p(R))$$
is zero.
\end{lem}

\begin{proof}
In the proof of Theorem 5.11 in \cite{Qui73}, we have the following formula
$$K_n(\MMM^{p+1}(R))=\underset{t:\substack{\text{regular}\\ \text{element}}}{\colim} K_n(\MMM^p(R/tR)).$$
Since $R$ is UFD by \cite{AB59}, for any regular element $t$ in $R$, we can write $t=p^{e_1}_1p^{e_2}_2\ldots p^{e_r}_r$ where $p_i$ are prime elements. By d\'evissage theorem in \cite{Qui73}, we have the following formula 
$$K_n(\MMM^p(R/tR)) \isoto K_n(\MMM^p(R/p_1p_2\ldots p_r R)).$$
\noindent
\textbf{Claim}\\
We have the following formula
$$K_n(\MMM^p(R/p_1p_2\ldots p_r R)) \isoto \underset{i=1}{\overset{r}{\bigoplus}} K_n(\MMM^p(R/p_i R)).$$
\begin{proof}[Proof of \textbf{Claim}]
We put $X=\Spec R/p_1p_2\ldots p_r R$ and $X_i=\Spec R/p_iR$. For any closed set $Y \subset X_?$, we put $\Perf^Y(X_?)$ the category of strictly perfect complexes which are acyclic on $X_?-Y$. We also put $\Perf^p(X_?):=\underset{\substack{Y \subset X_? \\ \Codim Y \geqq p}}{\cup} \Perf^Y(X_?)$. Then we have the following identities
\begin{eqnarray*}
K_n(\MMM^p(X)) & \underset{\text{I}}{\isoto} & \underset{\substack{Y \subset X \\ \Codim Y \geqq p}}{\colim} K'_n(Y) \underset{\text{II}}{\isoto} \underset{\substack{Y \subset X \\ \Codim Y \geqq p}}{\colim} K_n(X \ \text{on}\ Y)\\
& \isoto & \underset{\substack{Y \subset X \\ \Codim Y \geqq p}}{\colim} \underset{i=1}{\overset{r}{\bigoplus}} K_n(X_i \ \text{on} \ X_i \cap Y) \underset{\text{III}}{\isoto} \underset{i=1}{\overset{r}{\bigoplus}} \underset{\substack{Y \subset X_i \\ \Codim Y \geqq p}}{\colim} K_n(X_i \ \text{on} \ Y)\\
& \underset{\text{II}}{\isoto} & \underset{i=1}{\overset{r}{\bigoplus}} \underset{\substack{Y \subset X_i \\ \Codim Y \geqq p}}{\colim} K'_n(Y) \underset{\text{I}}{\isoto} \underset{i=1}{\overset{p}{\bigoplus}} K_n(\MMM^p(X_i))
\end{eqnarray*}   
where the isomorphisms I are proved by continuity \cite{Qui73}, \cite{TT90}, the isomorphisms II are proved by the Poincar\'e duality and comparing the following fibration sequences \cite{Qui73} and \cite{TT90}
$$K'(Y) \to K'(X_?) \to K'(X_? -Y) ,$$
$$K(X_? \ \text{on} \ Y) \to K(X_?) \to K(X_? -Y) $$
for any closed set $Y \subset X$ and to prove the isomorphism III, we are using the fact that all $X_i$ are equidimensional.
\end{proof}  
Therefore to prove Gersten's conjecture we shall only check that for any prime element $f$, the inclusion map $\MMM^p(R/fR) \to \MMM^p(R)$ induces the zero map
$$K_n(\MMM^p(R/fR)) \to K_n(\MMM^p(R)).$$   
Since $R/fR$ is regular, inductive argument implies that to prove Gersten's conjecture we shall only check that for any regular sequence $f_1,\ldots,f_{p+1}$ such that $(f_1,\ldots,f_{p+1})$ is prime ideal, the inclusion map $\MMM(R/(f_1,\ldots,f_{p+1})) \to \MMM^p(R)$ induces the zero map 
$$K_n(\MMM(R/(f_1,\ldots,f_{p+1})) \to K_n(\MMM^p(R)).$$
Since $R/(f_1,\ldots,f_{p+1})$ is regular, we have $K_n(\PPP(R/(f_1,\ldots,f_{p+1})) \isoto K_n(\MMM(R/(f_1,\ldots,f_{p+1}))$ by resolution theorem in \cite{Qui73}. Hence we get the result. 
\end{proof}

\noindent
Now \textbf{Lemma \ref{Qui ind}} implies the following assertion by famous Gersten-Sherman argument in \cite{Ger73}, \cite{She82} p.240, which is an application of the universal property for algebraic $K$-theory associated with semisimple exact categories \cite{She92} Corollary 5.2. From now on let $\FFF$ be the category of finite pointed connected CW-complexes and frequently using the notations in \cite{Moc07}.

\begin{lem}[Gersten-Sherman reduction argument] \label{GS red arg} $ $
\\
To prove the main theorem, we shall only check the following assertion:\\
For any $X \in \FFF$, any non-negative integer $p$, and any regular sequence $f_1,\ldots,f_{p+1}$ in $R$, the canonical map induced from the inclusion map $\PPP(R/(f_1,\ldots,f_{p+1})) \hookrightarrow \MMM^p(R)$,
$$\tilde{R}_0(\pi_1(X),\PPP(R/(f_1,\ldots,f_{p+1}))) \to \tilde{R}_0(\pi_1(X),\MMM^p(R)) \overset{\tilde{\Sh}}{\to} [X,(\mathbb{K}(\MMM^p(R)))_0]_{\ast}$$
is zero.
\end{lem}

\begin{proof}
We have the following commutative diagram for each $X \in \FFF$:
$$\xymatrix{
\tilde{R_0}(\pi_1(X),\PPP(R/(f_1,\ldots,f_{p+1}))) \ar[r] \ar[d]_{\tilde{\Sh}} & \tilde{R_0}(\pi_1(X),\MMM^p(R)) \ar[d]_{\tilde{\Sh}}\\
[X,{(\mathbb{K}(\PPP(R/(f_1,\ldots,f_{p+1}))))}_0]_{\ast} \ar[r] & [X,{(\mathbb{K}(\MMM^p(R)))}_0]_{\ast}
}$$
It is well-known that $\mathbb{K}(\MMM^p(R))$ is a $H$-space, $\PPP(R/(f_1,\ldots,f_{p+1}))$ is semi-simple and by the universal property \cite{She92} Corollary 5.2, we learn that we shall only prove the composition
$$\tilde{R_0}(\pi_1(X),\PPP(R/(f_1,\ldots,f_{p+1}))) \to \tilde{R_0}(\pi_1(X),\MMM^p(R)) \overset{\tilde{\Sh}}{\to} [X,{(\mathbb{K}(\MMM^p(R)))}_0]_{\ast}$$
is the zero map for any $X \in \FFF$.
\end{proof}

\noindent
Next we will define equivalence relations between morphisms in $\MMM^p(R)$ as follows:\\
For any $R$-modules $M$, $N$ in $\MMM^p(R)$, and morphisms $f, g:M \to N$, we will declare $f \sim g$.\\
Then $\MMM^p(R)$ is an exact category with equivalence relations satisfying the cogluing axiom in the sense of \cite{Moc07}. So we can define the Grothendieck group of lax $G$-representations in $\MMM^p(R)$. (For the precise definition, see \cite{Moc07} Definition 3.9).

\begin{lem}[Retraction principle] \label{Ret pri} $ $
\\
To prove main theorem, we shall only check the following assertion:\\
In the notation \textbf{Lemma \ref{GS red arg}}, the canonical map induced from the inclusion map 
$$\PPP(R/(f_1,\ldots,f_{p+1})) \hookrightarrow \MMM^p(R),$$
$$R_0(G,\PPP(R/(f_1,\ldots,f_{p+1}))) \to R_0^{lax}(G,\MMM^p(R))$$
is zero.
\end{lem}

\begin{proof}
We have the following commutative diagram for each $X \in \FFF$:
$$\xymatrix{
\tilde{R_0}(\pi_1(X),\PPP(R/(f_1,\ldots,f_{p+1}))) \ar[r] & \tilde{R_0}(\pi_1(X),\MMM^p(R)) \ar[r]^{\tilde{\Sh}} \ar[d] & [X,{(\mathbb{K}^{\edge}(\MMM^p(R)))}_0]_{\ast} \ar[d]^{\text{I}}\\
& \tilde{R_0^{\lax}}(\pi_1(X),\MMM^p(R)) \ar[r]_{\tilde{\Sh}^{\lax}} & [X,{(\mathbb{K}^{\lax,\edge}(\MMM^p(R)))}_0]_{\ast}
}$$
where the morphism I is a injection by retraction theorem 3.13 in \cite{Moc07}. Hence we get the result.
\end{proof}

\noindent
The following argument is one of a variant of weight argument of the Adams operations. (See \cite{Moc07} \S 1.) 

\begin{lem}[Weight changing argument] \label{Wei chan} $ $
\\
The assertion in \textbf{Lemma \ref{Ret pri}} is true. Therefore Gersten's conjecture is true.
\end{lem}

\begin{proof}
We put $B=R/(f_1,\ldots,f_p)$. Let $G$ be a group and $(X,\rho_X)$ be a representation in $\PPP(B/f_{p+1}B)$. Since $B/f_{p+1}B$ is local, $X$ is isomorphic to ${(B/f_{p+1}B)}^{\oplus m}$ for some $m$ as a $B/f_{p+1}B$-module. Then there is a short exact sequence
$$0 \to B^{\oplus m} \overset{f_{p+1}}{\to} B^{\oplus m} \overset{\pi}{\to} X \to 0.$$
For each $g \in G$, we have a lifting of $\rho_X(g)$, that is, a $R$-module homomorphism $\tilde{\rho}(g):B^{\oplus m} \to B^{\oplus m}$ such that $\tilde{\rho}(g) \mod f_{p+1}=\rho_X(g)$. Since $[B^{\oplus m} \overset{f_{p+1}}{\to} B^{\oplus m}]$ is a minimal resolution of $X$ as a $B$-module, (For the definition of a minimal resolution, see \cite{Ser00} p.84.) we can easily learn that $\tilde{\rho}(g)$ is an isomorphism as a $B$-modules by Nakayama's lemma. Therefore $\tilde{\rho}(g)$ is an isomorphism as a $R$-modules. Obviously assignment $\tilde{\rho}:G \to \Aut(B^{\oplus m})$ defines a lax representation $(B^{\oplus m},\tilde{\rho})$ in $\MMM^p(R)$ and we have a short exact sequence
$$(B^{\oplus m},\tilde{\rho}) \overset{f_{p+1}}{\to} (B^{\oplus m},\tilde{\rho}) \overset{\pi}{\to} (X,\rho_X)$$
in $\LAX(\ulG,\MMM^p(R))_{\strict}$. Notice that proving $f_{p+1}$ is a strict deformation, we need the assumption that $R$ is commutative!! So we have an identity $$[(X,\rho_X)]=[(B^{\oplus m}, \tilde{\rho})]-[(B^{\oplus m},\tilde{\rho})]=0$$
in $R_0^{\lax}(G,\MMM^p(R))$. Hence we get the result. 
\end{proof}

\section{Corollaries}

\noindent
In this section, we will discuss applications of \textbf{Theorem \ref{Main theorem}}. First we get the following absolute version of Bloch's formula.

\begin{cor} $ $
\\
For a regular noetherian scheme $X$, there is a canonical isomorphism
$$\Homo^p(X,\KKK_p)\isoto A^p(X)$$ 
where $\KKK_p$ is the Zariski sheaf on $X$ associated to the presheaf $U \mapsto K_n(U)$ and $A_p(X)$ is defined by the following formula
$$A^p(X):=\coker(\underset{x \in X_{p-1}}{\coprod} k(x)^{\times} \overset{\ord_x}{\to} \underset{x \in X_p}{\coprod} \mathbb{Z}).$$
Here $X_i$ is the set of points of codimension $i$ in $X$.
\end{cor}

\begin{proof}
Combining Propositions 5.8 and 5.14 and Remark 5.17 in \cite{Qui73} and \textbf{Theorem \ref{Main theorem}}, we can easily obtain the result.
\end{proof}

\noindent
Next we will cite the following well-known statement.

\begin{prop} \label{equi cond} $ $
\\
{\rm (c.f. \cite{Lev85} P.452, Proposition 1.1, \cite{Moc07} Proposition 1.2)} Let $A$ be a commutative regular local ring. Then the following statements are equivalent.\\ 
{\rm (i)} The maps $K_0(\MMM^p(A)) \to K_0(\MMM^{p-1}(A))$ are zero for $p=1,\cdots,\dim A$.\\
{\rm (ii)} $K_0(\MMM^p(A))$ is generated by cyclic modules $A/(f_1,\cdots,f_p)$ where $f_1,\cdots,f_p$ forms a regular sequence for $p=1,\cdots,\dim A$.\\
{\rm (iii)} $A_p(\Spec A)=0$ for any $p < \dim A$.
\end{prop}

\noindent
Therefore we get the following results.

\begin{cor}[Vanishing conjecture] $ $
\\
For any commutative regular local ring $R$ and any $p < \dim R$, we have $A_p(\Spec A)=0$.
\end{cor}

\begin{cor}[Generator conjecture] \label{gen con} $ $
\\
For any commutative regular local ring $R$, $K_0(\MMM^p(R))$ is generated by cyclic modules $R/(f_1,\cdots,f_p)$ where $f_1,\cdots,f_p$ forms a regular sequence for $p=1,\cdots,\dim R$.\\
\end{cor}


\begin{thebibliography}{99}
\bibitem[AB59]{AB59}
M. Auslander and D. Buchsbaum, {\it{Unique factorization in regular local rings}}, Proc. Nat. Acad. Sci. USA., \textbf{45} (1959), p.733-734.
\bibitem[Ger73]{Ger73} 
S. Gersten, {\it{Some exact sequences in the higher K-theory of rings}}, In Higher K-theories, Springer Lect. Notes Math. \textbf{341} (1973), p.211-243.
\bibitem[Lev85]{Lev85}
M. Levine, {\it{A $K$-theoretic approach to multiplicities}}, Math. Ann. \textbf{271} (1985), p.451-458.
\bibitem[Moc07]{Moc07}
S. Mochizuki, {\it{Gersten conjecture for commutative discrete valuation rings}}, available at http://www.math.uiuc.edu/K-theory/0819 (2007).
\bibitem[Qui73]{Qui73}
D. Quillen, {\it{Higher algebraic K-theory I}}, In Higher K-theories, Springer Lect. Notes Math. \textbf{341} (1973), p.85-147.
\bibitem[Ser00]{Ser00}
J. P. Serre, {\it{Local algebra}}, Springer monographs in Mathematics (2000)
\bibitem[She82]{She82}
C. Sherman, {\it{Group representations and algebraic $K$-theory}}, In Algebraic $K$-theory, Part I (Oberwolf ach,1980), Springer Lect. Notes Math \textbf{966} (1982), p.208-243.
\bibitem[She92]{She92}
C. Sherman, {\it{Group representations and algebraic $K$-theory:II}}, In Contemporary Math. Vol. \textbf{126} (1992), p.165-177.
\bibitem[TT90]{TT90}
R. W. Thomason, T. Trobaugh, {\it{Higher K-theory of schemes and of derived categories}}, In The Grothendieck Festscrift,Vol III, (1990), p.247-435. 
\end{thebibliography}
\end{document}